\documentclass{amsart} 
\usepackage{amssymb,amsmath,amsfonts,amstext} 
\newtheorem{theorem}{Theorem}
\newtheorem{lemma}[theorem]{Lemma}
\title{Observations concerning G\"{o}del's$\text{ }$1931} 
\author{Paola Cattabriga}
\address{University of Bologna --  Italy}
\email{co14099@iperbole.bologna.it}
\begin{document}
\begin{abstract}
This article demonstrates  the invalidity of Theorem VI  in G\"{o}del's
monograph of 1931, by showing that 
\[
\begin{array}[t]{ll}
(15)
& \quad \quad \overline{ x B_{\kappa}(17\, \text{Gen} \, r)} \longrightarrow   
\text{Bew}_{\kappa}[Sb(r^{17}_{Z(x)})],  \\
 (16) &\quad \quad x B_{\kappa}(17\, \text{Gen} \, r)  \longrightarrow   
\text{Bew}_{\kappa}[\text{Neg}(Sb(r^{17}_{Z(x)})],
\end{array}\] 
(derived by 
definition (8.1)
$Q(x,y)\equiv \overline{x B_{\kappa}[Sb(y^{19}_{Z(y)})]}$ respectively from 
\[
\begin{array}[t]{ll}
(3) & R(x_{1},\ldots , x_{n}) \longrightarrow   
\text{Bew}_{\kappa}[Sb(r^{u_{1}\ldots u_{n}}_{Z(x_{1})\ldots Z(x_{n})})],  \\
(4) & \overline{R}(x_{1},\ldots , x_{n})  \longrightarrow   
\text{Bew}_{\kappa}[\text{Neg}(Sb(r^{u_{1}\ldots u_{n}}_{Z(x_{1})\ldots Z(x_{n})})],
\end{array}\]
of Theorem V) 
are false in $P$. This is achieved in two steps.  
First, the predicate complementary to the well-known G\"{o}del's predicate $\text{Bew}(x)$ is
defined by adding a new relation $\text{Wid}(x)$. In accordance, new logical
connections are established, Lemma (6).
Second,
\[
\begin{array}[t]{ll}
(I)  & \quad \quad
 \overline{\overline{ x B_{k}(17 Gen r)} \rightarrow 
Bew_{k}[Sb(r^{17}_{Z(x)})]},\\
(II) & \quad \quad \overline{x B_{\kappa}(17 Gen r)\rightarrow 
Bew_{k}[Neg(Sb(r^{17}_{Z(x)})]}, 
\end{array}\]
are derived from Lemma (6) and definition (8.1). It amounts to saying that (15) and (16) are false and unacceptable
for the system ($P$).
On the  account of that, the two well-known cases
\[\begin{array}[t]{ll}
1. \quad 17\, \text{Gen} \, r \text{ is not } \kappa \text{--PROVABLE}, 
\\ 
 2. \quad \text{Neg}(17\, \text{Gen} \, r) \text{ is not } \kappa \text{--PROVABLE},
\end{array}\]
can not be drawn (unless we  say that they are true 
because 
$\, ( (15) \, \& \, (I)) \supset 1.  \, $  and  $ \, ( (16) \, \& \, 
(II) ) \supset 2. \, $ 
are examples of $(p \, \&\,  \text{\tiny{$\sim$}}  p)\supset q$, i. e. 
{\it ex falso sequitur quodlibet}). 
\end{abstract}
    \maketitle

\subsection*{Introduction}\label{intro}
``Wir m\"{u}ssen wissen. Wir werden wissen." We must know. We will know.
With these words, pronounced at the  International Congress of 
Mathematicians in Paris, in 1900, Hilbert settled the basis of the logical 
investigations of the foundations of mathematics of the twentieth century. 
Whatever they were commands or wishes, they put 
the human spirit in a position of certainty about the solvability of 
every mathematical problem. 
These words embodied  ancient aspirations, such as the achievement of a general problem-solving method, 
like  Ramon Lull's \emph{Ars Magna}, or the completion of a universal symbolic
language like Leibnitz's \emph{Characteristica} and Frege's \emph{Begriffsschrift} 
\cite[-4]{borger}.
 But for  histories sense of humor, or 
for the unpredictability of  human nature, they were just at the 
starting point of the exploration which led to the accomplishment
of the opposite position, see for example the conviction 
of Godfrey Hardy in 1928.
\begin{quote}  
`` There is of course no such theorem, and this is very fortunate, since 
if there were we should have a mechanical set of rules for the 
solution of all mathematical problems, and our activities as 
mathematicians would come to an end."   \cite[-93]{hodges} 
\end{quote} 
It was hard to believe that statements about numbers like Fermat's Last Theorem and 
Goldbach's conjecture which the efforts 
of the centuries have failed  to solve,  could in fact be decided  by a mechanical 
process. All this  found  its culmination in  G\"{o}del's 1931 
article, that placed severe limits on the power of 
mathematical reasoning and on the power of the axiomatic method \cite{godel1}.

Many mathematicians at the end of the nineteenth-century considered  
consistency to be
sufficient in securing the existence of theoretical objects. The result 
obtained by
 G\"{o}del in 1930  assuring the existence of a model for the first 
order logic,  apparently   confirmed  this opinion \cite{godelc}. 
However,   for 
G\"{o}del    
comparing  consistency with  existence  
\begin{quote}  
`` manifestly presupposes the axiom that every mathematical problem is 
solvable. Or, more precisely, it presupposes that we cannot prove the 
unsolvability of any problem." \cite[60-61]{godelt}   
\end{quote}
Through 
 G\"{o}del's contributions, the  reference to the existence of unsolvable problems in 
 mathematics is recurring (see note 61 in \cite[190-191]{godel1}), together with it's winding 
exchange  with the existence of undecidable 
propositions \cite[144-145]{godel1}.
This connection of unsolvable problems  and 
undecidable propositions, was supported by the 
conviction that  
\begin{quote}  
`` there are true propositions (which may even 
be provable by means of other principles) that cannot be derived in 
the system under consideration." \cite[102-103]{godelc} 
\end{quote}
    So that, for G\"{o}del, there 
were unsolvable problems in mathematics that, although contentually true,   
 were unprovable in the formal system  
 ($P$ in \cite[150-151]{godel1}). 
\begin{quote}  
``(Assuming the consistency of classical mathematics) one can even 
give examples of propositions (and in fact of those of the type of 
Goldbach or Fermat) that, while contentually true, are unprovable in 
the formal system of classical mathematics." \cite[202-203]{godel31a} 
 
\end{quote}
Nowadays the developments concerning Fermat's last theorem 
\cite{wiles1,wiles2},
clarify that \emph{problems with no solution} are  not so 
\emph{unprovable} in   
mathematics. We are then allowed to state that, in 
G\"{o}del's convictions, around the existence of undecidable problems 
in mathematics, there is a presupposition of unprovability that the 
present article will help to clarify.

Of all the remarkable logical achievements of the twentieth 
century perhaps the most outstanding is    G\"{o}del's celebrated
incompleteness argumentation of 1931.
In contrast to Hilbert's program which called for embodying classical
mathematics in a formal system and proving that system
consistent by finitary methods \cite{hilbert1,hilbert2},
G\"{o}del's paper showed that not even the first step could be carried 
out fully,  any formal system suitable for the arithmetic of 
integers was incomplete \cite{godel1,godel2}.

G\"{o}del's incompleteness argument holds today 
the same scientific status as Einstein's principle of relativity, 
Heisenberg's uncertainty principle, and Watson and Crick's double helix 
model of DNA \cite[-315]{calude}.
With 
absolute respect and acknowledgement for the extraordinary contribution given by  
G\"{o}del to the logical investigation, this article 
brings  G\"{o}del's achievement into question  by the definition 
of the refutability predicate as a number-theoretical statement.

As is well known, G\"{o}del numbering is a special device  by  which each expression in 
arithmetic can refer to itself.
Just a declarative programmer can catch the extraordinary 
expressive power of G\"{o}del numbering.
\begin{quote} 
    ``Many of the logicians who worked in the first half of this 
    century were actually the first inventors of programming 
    languages. G\"{o}del uses a scheme for expressing algorithms that 
    is much like the high-level language LISP ... ."  
 \cite[-19]{chaitin} 
``But if you read his original 1931 paper with the benefit of 
hindsight, you'll see that G\"{o}del is programming in LISP, ... ."  
\cite[-14]{chaitin} 
\end{quote}
It was precisely  with  this view of G\"{o}del   
 programming recursive statements that the content of this work  
was conceived.
According to  G\"{o}del's first 
incompleteness argument it is possible to construct a formally 
undecidable proposition in PM,  a
statement that,  although true,  turns out to be neither provable  nor 
refutable for the system. 
This article develops proof of  invalidity 
of Theorem VI in G\"{o}del's 1931,
the so-called G\"{o}del's  first incompleteness theorem, in two steps: defining 
 refutability within the same recursive status as provability  and 
showing  that as a consequence propositions (15) and (16), derived from  
definition 8.1 in G\"{o}del's 1931  \cite[174-175]{godel1}, are 
false and  unacceptable for the system. The achievement of their falsity 
blocks the derivation of Theorem VI, which turns out to be  therefore invalid,
together with all the depending  theorems  (Sections \ref{wider},  
 \ref{antig}  and \ref{finalsynt}).

Completeness is connected to decidability  since a resolution 
procedure takes advantages of those   sorts of logical inferences   
 formalized also by the axiomatic systems. 
With respect to the finitely 
axiomatized  systems, it leads  to a theoretical list of all  
theorems, which is  in turn  a procedure of semi-decidability  for 
the first order logic.
With regard to that, the  Lemmas  of Section \ref{wider}     are clarifying  and  
opening new perspectives,   partially already exposed in 
\cite{catta1,catta2}.
In fact   this work  applies  directly to the text of 
G\"{o}del's 1931 article, the results obtained in 
\cite{catta1} for the modern 
version of the incompleteness argument, which is mainly based on the 
diagonalization lemma. It provides therefore   
as a novel contribution   
clear and advanced proof of the invalidity of   G\"{o}del's   original
 argument.
 For all these aspects the result  exposed here stands out from previous 
 criticism of G\"{o}del's 1931, mainly focused on the antinomic 
 features of G\"{o}del's self-reference statement, (see  
 \cite{kucz}, \cite{barzin},  \cite{perelman}, 
 resumed in  \cite{ladriere};  
 \cite{rivetti}).
In the following  the reader is required to be familiar with   
 G\"{o}del's 1931  \cite{godel1}.
 \section{}\label{wider}
 Let us begin adding to the list of functions (relations) 1-45 in 
G\"{o}del's 1931 two new relations, 45.1 and 46.1, in terms of the 
preceding ones by the procedures given in Theorems I-IV \cite[158-163]{godel1}. 
We shall recall only the well-known definitions 44., 45. and 
46., for the whole list the reader is referred to \cite[162-171]{godel1}.

\begin{math}
\text{44. }Bw(x)\equiv  
     (n)\{ 0<n\leq l(x)  \longrightarrow Ax(n\, Gl\, x) \vee        
     (Ep,q)[0<p,q<n \;\&\;  \\
        \quad \qquad \text{   }  \quad \text{   } \qquad \text{   }\qquad \text{   }  \quad 
   Fl(n\, Gl\, x, p \,Gl\, x, q \,Gl\, x )] \} \;\&\;  l(x) > 0,  \\
\text{$x$ is a PROOF ARRAY } 
\end{math}
(a finite sequence of FORMULAS, each of which is 
either an AXIOM or an IMMEDIATE CONSEQUENCE of two of the preceding 
FORMULAS).

\begin{math}  
\text{45. } xBy \equiv Bw(x)\;\&\; [l(x)]\, Gl\, x = y, \\
\text{$x$ is a PROOF of the FORMULA y.}
\end{math}

\begin{math}  
\text{45.1. } xWy \equiv Bw(x)\;\&\; [l(x)]\, Gl\, x = \text{Neg}(y),\\
\text{$x$ is a REFUTATION of the FORMULA y.}
\end{math}

\begin{math}  
\text{46. } \text{Bew}(x) \equiv (Ey) y B x ,\\
\text{$x$ is a PROVABLE FORMULA. }
\end{math}
($\text{Bew}(x)$ is 
the first one of the notions 
1--46 of which we cannot assert that it is recursive.)

\begin{math}  
\text{46.1. } \text{Wid}(x) \equiv(Ez) z W x,\\
\text{$x$ is a REFUTABLE FORMULA. }
\end{math}
($\text{Wid}(x)$ is the second one of the notions 
1--46 of which we cannot assert that it is recursive.) $\text{Wid}$ 
is the shortening for ``Widerlegung" and must not be 
mistaken with the notion defined by G\"{o}del in note 63 referring instead to 
``Widerspruchsfrei" \cite[192-193]{godel1}, which afterwards  
we will call   $\text{Wid}_{s}$.

Being classes included among relations, as one-place relation, 44. 
define the recursive class of the proof arrays. 
Recursive relations $R$ have the property that for every given $n$-tuple of numbers it 
can be decided whether $R(x_{1}\ldots x_{n})$ holds or not.
45. and 
45.1 define  recursive relations, $xBy$ and $xWy$, so that for every 
given couple of numbers it can be decided whether $xBy$ and $xWy$ hold 
or not. 
In accordance we can state the following Lemmas.

\begin{lemma}\label{1*}
 
$$
    (x)\{\text{\rm Wid}(x) \;\thicksim\;\text{\rm  Bew}[\text{\rm 
    Neg}(x)]\}  
\text{ and } 
     (x)\{\text{\rm Bew}(x)\; \thicksim\;\text{\rm Wid}[\text{\rm Neg}(x)]\}.
$$

\end{lemma}

\begin{proof}
 Immediately by  46.1., 46., 45.1., 45. and $\text{\rm  
 Neg}(\text{\rm  Neg}(x)) \thicksim x$.
 \end{proof}

 \bigskip 
 
 \begin{lemma}\label{notboth1}
For any $x$, $y$ in $P$, $$ \text{ not both } xWy \text{  and } xBy.$$ 
\end{lemma}
\begin{proof}
Let us suppose to have, for an arbitrary couple $x$, $y$, both $xWy$  and  $xBy$ in $P$. By 45. and 45.1 
\begin{equation*}
 Bw(x)\;\&\;[l(x)]\, Gl \, x = y  \quad \&  \quad 
 Bw(x)\;\&\;[l(x)]\, Gl \, x = \text{Neg}(y),
\end{equation*}
then by 44.
\begin{multline*}
 (n)\{  0<n\leq  l(x) \longrightarrow  Ax(n\, Gl\, x) \vee \\
(Ep,q)[0<p, q<n \;\&\;     Fl(n\, Gl\, x,  p \,Gl\, x, q \,Gl\, x )] \} \;\&\;  \\
l(x) > 0 \; \& 
\; [l(x)]\, Gl \, x = y \; \&   \;
[l(x)]\, Gl \, x = \text{Neg}(y).
\end{multline*}
For $l(x) > 0$ and $n= l(x)$   this should be
\begin{multline*}
 Ax(y) \;\&\; Ax(\text{Neg}(y)) \vee  
  (Ep,q)[0<p,q<n \;\&\;\\
  Fl(y, p \,Gl\, x, q \,Gl\, x )\;\&\; 
 Fl(\text{Neg}(y), p \,Gl\, x, q \,Gl\, x )]  
\end{multline*}
which is impossible by  42. $Ax(x)$, 43. $Fl(x,y,z)$ 
and  by the definitions of the axioms and of immediate consequence in $P$.
 \end{proof}
Recursive functions have the property that, for each given set of 
values of the arguments, the value of the function can be computed by 
a finite procedure \cite[-348, 369-371]{godel2}. The following Lemma, 
especially as regard to $l$, $Gl$  and $\text{Neg}$, 
is  involved with this property  and 
the today's notion of effective computability.

\bigskip 

\begin{lemma}\label{aut0}
For any $x, y$ such that  $\text{$x$ is a PROOF of the FORMULA y}$, in $P$, $$xBy \text{ or } xWy.$$
\end{lemma}
\begin{proof}
Assume that  $\text{$x$ is a PROOF of the FORMULA y}$. 
Then, 45., $xBy$, and from the  schema
$p \supset p \vee q$, II.2. \cite[154-155]{godel1}, 
 $xBy \supset xBy \vee  xWy$ so that  $ xBy \vee  xWy$,  
i. e. for any $x, y$ such that $\text{$x$ is a PROOF of the FORMULA y}$,  $ xBy \vee  xWy$.
 \end{proof} 
In the case that    $\text{$x$ is a PROOF of the FORMULA Neg(y)}$, 
i.e. $[l(x)]\, Gl\, x = \text{Neg}(y)$, the same 
Lemma can be derived as follows. Given $Bw(x)$,  
$[l(x)]\, Gl\, x$ is a natural number for a FORMULA $z$. $\text{Neg}$ is a 
recursive function so that we have a finite procedure to determine the 
value of $\text{Neg}(y)$. Let us suppose that $z = \text{Neg}(y)$, then 
 $xBz$ is 
true, and from $xBz \supset xBz  \vee  xWz$ we draw  $ xBz  \vee   
xWz$. We notice that this Lemma holds 
for any FORMULA $z$ such that $[l(x)]\, Gl\, x  = z$
even if $xBz$ is $xWy$.

\bigskip 

\begin{lemma}\label{wb}
For any $x, y$  such that $\text{$x$ is a PROOF of the FORMULA y}$,  in $P$,
$$xWy \thicksim   \overline{xBy}.$$
\end{lemma}
\begin{proof}
 Immediately by Lemma (\ref{notboth1}) and Lemma (\ref{aut0}).
\end{proof}
  Let us notice that  Lemma (\ref{wb}) yields, for any $x, y$ such that 
  $\text{$x$ is a PROOF }$   $\text{of the FORMULA y}$
  in $P$, 
also  $\overline{ xWy }\thicksim  xBy$.

\bigskip 

\begin{lemma}\label{ifneg}
For any $x, y$ such that $\text{$x$ is a PROOF of the FORMULA y}$, in $P$,
$$\overline{\text{Bew}(y)} \text{ if } \overline{ xBy} \text{ and }  \overline{\text{Wid}(y)} \text{ iff }  xBy.$$
\end{lemma}
\begin{proof}
Let us  assume that \text{$x$ is a PROOF of the FORMULA y} and 
$\overline{xBy}$.
As an immediate 
consequence we have $(x)\overline{xBy}$.
 By  46. $\text{Bew}(y) \thicksim (Ex)xBy$,  $\overline{\text{Bew}(y)} \thicksim 
\overline{(Ex)xBy}$ and  $\overline{\text{Bew}(y)} \thicksim 
(x)\overline{xBy}$. 
Accordingly
  for any $x, y$ such that 
 $\text{$x$ is a PROOF }$ $\text{of the FORMULA y}$, 
 $\overline{\text{Bew}(y)} \text{ if } \overline{xBy}.$  

Let us assume that \text{$x$ is a PROOF of the FORMULA y}. 
By 45. $xBy$ and from Lemma (\ref{wb})
$\overline{ xWy }\thicksim  xBy$, so that  $\overline{ xWy }$.   
$(x)\overline{ xWy }$ is an immediate consequence of $\overline{ 
xWy }$, hence  $(x)\overline{ xWy }$.  $(p \, \&\, q) 
\supset  ( p  \thicksim  q)$ yields then
 $(x)\overline{ xWy } \thicksim  xBy$.
By 46.1  $\text{Wid}(y) \thicksim$ $ (Ex)xWy$, so that  
 $\overline{\text{Wid}(y)}$ $\thicksim 
\overline{(Ex)xWy}$ and  $\overline{\text{Wid}(y)} \thicksim (x)\overline{xWy}.$
 Accordingly, for any $x, y$ such that 
 $\text{$x$ is a PROOF of the FORMULA y}$,
 $\overline{\text{Wid}(y)} \thicksim xBy.$
\end{proof}

All preceding Lemmas were carried out constructively, needlessly to assume
consistency.
Let us recall the following G\"{o}del's definitions \cite[-173]{godel1}.
Let $\kappa$ be any class of FORMULAS and $Flg(\kappa)$  the 
smallest set of FORMULAS that contains all FORMULAS of $\kappa$ and 
all AXIOMS, and is closed under the relation IMMEDIATE CONSEQUENCE.
\begin{equation}\tag{5}\label{cinque}
\begin{split} 
     Bw_{\kappa}(x)\equiv &(n)[n \le l(x) \longrightarrow   Ax(n\, Gl \, x) 
  \vee   (n\, Gl \, x) \in \kappa \vee \\ &
   (Ep,q) \{ 0<p,q<n \;\&\;  
    Fl(n\, Gl x, p\, Gl x, q Gl x)\} ]   \;\&\; l(x) > 0, 
\end{split}
\end{equation}
\begin{equation}\tag{6}\label{sei}
    xB_{\kappa}y\equiv Bw_{\kappa}(x)\;\& \;[l(x)]\, Gl\, x = y,
\end{equation}
\begin{equation}\tag{6.1}\label{6.1}
\text{Bew}_{\kappa}(x)\equiv (Ey)y B_{\kappa} x,
\end{equation}
\begin{equation}\tag{7}\label{7}
(x)[\text{Bew}_{\kappa}(x) \thicksim x \in Flg(\kappa)].
\end{equation}
Let us  augment this list with two new definitions.
\begin{equation}\tag{6.2}\label{62}
xW_{\kappa}y \equiv Bw_{\kappa} (x)\; \& \;[l(x)] \, Gl \, x = \text{Neg}(y),
\end{equation}
\begin{equation}\tag{6.2.1}\label{621}
\text{Wid}_{\kappa}(x)\equiv (Ey)y W_{\kappa} x.
\end{equation}
We add further that  (\ref{6.1}) and (\ref{7}) yield
\begin{equation}\tag{7.1}\label{7.1}
(x)[ (Ey)y B_{\kappa} x \thicksim x \in Flg(\kappa)].
\end{equation}
\bigskip 

\begin{lemma}\label{ifnegK}
For any $x$    such that $Bw_{\kappa}(x)$ and any  FORMULA $ y$ such 
that $\quad $ 
 $[l(x)]\, Gl\, x = y$ in $\kappa$
 
 $$\overline{\text{Bew}_{\kappa}(y)} \text{ if }\overline{ xB_{\kappa}y}  
 \text{ and } 
 \overline{\text{Wid}_{\kappa}(y)} \text{ iff } xB_{\kappa}y.
$$
\end{lemma}
\begin{proof}
 By Lemma (\ref{ifneg}) and previous definitions.
\end{proof}

\bigskip 

\begin{lemma}\label{infersub}
Given a  CLASS SIGN $a$ with the FREE VARIABLE $v$ and
 a SIGN $c$ of the same type as $v$,
    \begin{gather} 
   \text{  if } \quad    
    \overline{\text{Bew}_{\kappa}(v\, \text{Gen}  \, a)} \quad  
    \text{ then } \quad    
    \overline{\text{Bew}_{\kappa}[Sb(a^{v}_{c})]},\tag{i} \\
\tag{ii}
  \text{ if } \quad  
  \overline{\text{Wid}_{\kappa}(v\, \text{Gen}  \, a)} \quad  
    \text{ then }  \quad    
    \overline{\text{Wid}_{\kappa}[Sb(a^{v}_{c})]}. \quad  
\end{gather}
\end{lemma}
\begin{proof}

(i) Let us suppose that $\text{Bew}_{\kappa}[Sb(a^{v}_{c})],$ then  
from 
 (\ref{7}) $Sb(a^{v}_{c})  \in Flg(\kappa).$ 
$Flg(\kappa)$ is closed under the relation IMMEDIATE CONSEQUENCE 
therefore
  $ v\, \text{Gen}  \, a \in Flg(\kappa).$ Thus, by (\ref{7}) 
  again, 
$\text{Bew}_{\kappa}(v\, \text{Gen}  \, a).$ Accordingly  if 
$\overline{\text{Bew}_{\kappa}(v\, \text{Gen}  \, a)}$   then     
    $\overline{\text{Bew}_{\kappa}[Sb(a^{v}_{c})]}$.

(ii) Let us suppose that $\overline{\text{Wid}_{\kappa}(v\, \text{Gen}  
\, a)}$, then, by Lemma (\ref{ifnegK}), for any $x$ such that $Bw_{\kappa}(x)$  and 
 $[l(x)]\, Gl\, x =   v\, \text{Gen}  \, a $ in $\kappa$,
 $x B_{\kappa} (v\, \text{Gen}  \, a) $. By 
 (\ref{7.1})\footnote{
 If, for any $x$ such that $Bw_{\kappa}(x)$  and 
 $[l(x)]\, Gl\, x =   v\, \text{Gen}  \, a $, 
 $x B_{\kappa} (v\, \text{Gen}  \, a) $ then, for any $x$ such that $Bw_{\kappa}(x)$  and 
 $[l(x)]\, Gl\, x =   v\, \text{Gen}  \, a $, 
 $(Ey)y B_{\kappa} (v\, \text{Gen}  \, a) $.
 } we have 
 therefore
  $v\, \text{Gen}  \, a \in Flg(\kappa)$. Thanks to Axiom 
III.1\footnote{
  \text{III. } \text{ Any formula that results from the schema}   
$$ \text{1. }  v \Pi (a) \supset \text{Subst }  a(^{v}_{c})   $$
 when the following substitutions are made for $a$, $v$, 
    and $c$ (and the operation indicated by  Subs is performed in 
    1):
for $a$ any formula, for $v$ any variable and for $c$ any 
    sign of the same type as $v$, provided $c$
   does not contain any 
    variable that is bound in $a$ at a place where $v$ is free
\cite[154-155; 37. 38. 168-169]{godel1}. 
},
we obtain   $  Sb(a^{v}_{c}) \in Flg(\kappa)$ and
from (\ref{7.1}) again, we have $(Ez)z B_{\kappa} [Sb(a^{v}_{c})]$. 
Finally, 
 from Lemma (\ref{ifnegK}) we have
$ \overline{\text{Wid}_{\kappa}[Sb(a^{v}_{c})]} $
(for any $x$ such that $Bw_{\kappa}(x)$  and  $[l(x)]\, Gl\, x =   v\, \text{Gen}  
  \, a $).
  \end{proof}

  \section{}\label{antig}
We are now ready to derive the main result of this article, Theorem 
(\ref{no-undec}). 
The invalidity of Theorem VI in G\"{o}del's  1931 article  
\cite[172-177]{godel1} 
follows from Lemmas (\ref{ifnegK}) and (\ref{infersub}).
In G\"{o}del's  1931 argumentation the proof that both $17\, \text{Gen} \, r$ and 
$\text{Neg}(17\, \text{Gen} \, r)$ are not 
$\kappa\text{--PROVABLE}$ is based on the two statements \cite[174-176]{godel1}
\begin{equation}\tag{15}\label{15}
 \overline{ x B_{\kappa}(17\, \text{Gen} \, r)} 
 \longrightarrow   \text{Bew}_{\kappa}[Sb(r^{17}_{Z(x)})]
	\end{equation}
 \begin{equation}\tag{16}\label{16}
x B_{\kappa}(17\, \text{Gen} \, r) 
 \longrightarrow   \text{Bew}_{\kappa}[\text{Neg}(Sb(r^{17}_{Z(x)})],	
\end{equation}
 which are respectively deduced from $Q(x,y)$ and 
$\overline{Q}(x,y)$, whereas
\begin{equation}\tag{8.1}\label{8.1}
Q(x,y)\equiv \overline{x B_{\kappa}[Sb(y^{19}_{Z(y)})]}.
\end{equation}
More precisely, $Q(x,y)$ is an instance of 
$R(x_{1},\ldots,x_{n})$ in 
 \begin{equation}\tag{3}
 R(x_{1},\ldots , x_{n}) \longrightarrow   
\text{Bew}_{\kappa}[Sb(r^{u_{1}\ldots u_{n}}_{Z(x_{1})\ldots Z(x_{n})})],  
\end{equation}
as exactly as $\overline{Q}(x,y)$ is an example of
$\overline{R}(x_{1},\ldots,x_{n})$ in 
 \begin{equation}\tag{4}
 \overline{R}(x_{1},\ldots , x_{n})  \longrightarrow   
\text{Bew}_{\kappa}[\text{Neg}(Sb(r^{u_{1}\ldots u_{n}}_{Z(x_{1})\ldots Z(x_{n})})],
\end{equation}
(Theorem V \cite[170-171]{godel1}),
so that  (3) $\longrightarrow$ (9) $\longrightarrow$ (15) and 
(4) $\longrightarrow$ (10) $\longrightarrow$ (16) \cite[170-175]{godel1}.

In accordance with  Lemmas (\ref{ifnegK}) and (\ref{infersub}) we will 
show that   given 
$Q(x,y)$,    (\ref{15}) turns out to be false,  and, similarly, 
given $\overline{Q}(x,y)$,  (\ref{16}) results to be false. 
We will show then as a consequence that Theorem VI in G\"{o}del's  
1931 \cite[172-173]{godel1}  is not achievable.

\bigskip 

\proof
 Let us assume $Q(x,y)$ in $\kappa$, then by definition (8.1)
 \begin{equation}\tag{I.1}\label{stepa}
 \overline{x B_{\kappa}[Sb(y^{19}_{Z(y)})]}.
 \end{equation}
 We substitute in it $p$ for y, see definitions 
 \begin{equation}\tag{11}
 p = 17\, \text{Gen}\, q, 
 \end{equation}
 ($p$ is a CLASS SIGN with the FREE VARIABLE 19),
 \begin{equation}\tag{12}
 r = Sb(q^{19}_{Z(p)}) 
  \end{equation}
  ($r$ is a recursive CLASS SIGN with the FREE VARABLE 17)
 and
  \begin{equation}\tag{13}
Sb(p^{19}_{Z(p)}) = 
Sb([17\, \text{Gen} \, q]^{19}_{Z(p)}) =  17 \, \text{Gen} \, Sb(q^{19}_{Z(p)}) = 17 
\, \text{Gen} \, r 
 \end{equation}
in  \cite[174-175]{godel1}, so that we have 
\begin{equation}\tag{I.2}\label{step1}
 \overline{ x B_{\kappa}(17\, \text{Gen}  \, r)}.
 \end{equation}
 By $a \rightarrow (\, _{ \widetilde{ \: } } b \rightarrow \,  _{ \widetilde{ \: } } (a \rightarrow b)) $  
 we have
 \begin{multline*}\tag{I.3}
 \overline{x B_{\kappa}(17\, \text{Gen}  \, r)} \longrightarrow 
 \Big[ \; \overline{\text{Bew}_{\kappa}[Sb(r^{17}_{Z(x)})]} 
 \longrightarrow  \\
 \overline{ 
 \overline{x B_{\kappa}(17\, \text{Gen}  \, r)} \longrightarrow 
 \text{Bew}_{\kappa}[Sb(r^{17}_{Z(x)})\;  \big ] }\;\Big ],
 \end{multline*}
 hence, by (I.2),
\begin{equation}\tag{I.4}
 \overline{\text{Bew}_{\kappa}[Sb(r^{17}_{Z(x)})]}  \longrightarrow  
 \overline{\; 
 \overline{x B_{\kappa}(17\, \text{Gen}  \, r)} \longrightarrow 
 \text{Bew}_{\kappa}[Sb(r^{17}_{Z(x)})]\; }.
\end{equation}
 Lemma (\ref{ifnegK}) yields 
 for any $x$ such that $Bw_{\kappa}(x)$, and $[l(x)]\, Gl\, x = $ $17\, \text{Gen}  \, 
 r$, 
\[\overline{\text{Bew}_{\kappa}(17\, \text{Gen}  \, r)}  \text{ if }
\overline{xB_{\kappa}(17\, \text{Gen}  \, r)},\] 
 accordingly from (\ref{step1}) 
\begin{equation}\tag{I.5}
    \overline{\text{Bew}_{\kappa}(17\, \text{Gen}  \, r)},
    \end{equation}
and by Lemma (\ref{infersub}) (i) 
\begin{equation}\tag{I.6}
    \overline{\text{Bew}_{\kappa}[Sb(r^{17}_{Z(x)})]}.
    \end{equation}   
(I.4) and (I.6) yield, for any $x$ such that $Bw_{\kappa}(x)$, and $[l(x)]\, Gl\, x = 17\, \text{Gen}  \, 
 r$,
 \begin{equation}\label{anti1}\tag{I}
 \overline{ \overline{ x B_{\kappa}(17\, \text{Gen} \, r)} 
 \longrightarrow   \text{Bew}_{\kappa}[Sb(r^{17}_{Z(x)})]}.
\end{equation}
\qed

\bigskip 

\begin{proof}
Let us assume $\overline{Q}(x,y)$ in $\kappa$, so that by (8.1), substituting $p$ 
for $y$, (11), (12) and (13), we obtain
\begin{equation}\tag{II.1}\label{step2}
x B_{\kappa}(17\, \text{Gen} \, r) .
\end{equation}
From $a \rightarrow 
 (\text{\tiny{$\sim$}} b \rightarrow \text{\tiny{$\sim$}} (a 
 \rightarrow b)) $,  
\begin{multline*}\tag{II.2}
x B_{\kappa}(17\, \text{Gen} \, r)  \longrightarrow \\
 \Big[\;\overline{\text{Wid}_{\kappa}[Sb(r^{17}_{Z(x)})]} \longrightarrow  
 \overline{ \; x B_{\kappa}(17\, \text{Gen} \, r) 
 \longrightarrow  \text{Wid}_{\kappa}[Sb(r^{17}_{Z(x)})]\;  }\Big ],
\end{multline*}
and then
\begin{equation}\tag{II.3}
\overline{\text{Wid}_{\kappa}[Sb(r^{17}_{Z(x)})]}   \longrightarrow  
 \overline{ \; x B_{\kappa}(17\, \text{Gen} \, r) 
 \longrightarrow  \text{Wid}_{\kappa}[Sb(r^{17}_{Z(x)})]\; }.
\end{equation}
From Lemma (\ref{ifnegK}), for any $x$ such that $Bw_{\kappa}(x)$, and $[l(x)]\, 
Gl\, x = 17\, \text{Gen}  \, r$ $ \text{ }$
\[\overline{\text{Wid}_{\kappa}(17\, \text{Gen} \, r)} \text{ iff }
x B_{\kappa}(17\, \text{Gen} \, r),\] and II.1
\begin{equation}\tag{II.4}
\overline{\text{Wid}_{\kappa}(17\, \text{Gen} \, r)}.
\end{equation}

Lemma (\ref{infersub}) 
(ii)\footnote{
$Z(x)$   does not contain any  variable bound in $r$ at a place at which $17$ is free. This mainly because 
$Z(x)$ is by definition,  $\text{ 17. } Z(n)$  \cite[164-165]{godel1}, a  NUMERAL and as such 
does not contain any variable. But if we think to $x$ in $Z(x)$ as to the 
PROOF ARRAY such that
$Bw_{\kappa}(x)$ and $[l(x)]\, Gl\, x = $ $17\, \text{Gen}  \,  r$, 
$x$ cannot be a variable bound in $r$ at a place at which $17$ is free, since
 $x  \, > \, r$ (see definitions $ \text{6. } n \, Gl\, x$    
 and $\text{ 7. } l(x)$  \cite[162-163]{godel1}) and, as it can be shown  by induction on the 
 recursive definition  $\text{ 44. } Bw(x)$, $x$ is 
 not  in $r$.
} yields, for any $x$ such that $Bw_{\kappa}(x)$, and $[l(x)]\, 
Gl\, x = 17\, \text{Gen}  \, r$,
\begin{equation}\tag{II.5}
\overline{\text{Wid}_{\kappa}[Sb(r^{17}_{Z(x)})]}
\end{equation}
and from II.3 
\begin{equation}\tag{II.6}
\overline{x B_{\kappa}(17\, \text{Gen} \, r) 
 \longrightarrow   \text{Wid}_{\kappa}[Sb(r^{17}_{Z(x)})]}.
\end{equation}
Finally, by Lemma (\ref{1*}),  for any $x$ such that 
$Bw_{\kappa}(x)$, and $[l(x)]\, Gl\, x = 17\, \text{Gen}  \, r$,
\begin{equation}\label{anti2}\tag{II}
	\overline{  x B_{\kappa}(17\, \text{Gen} \, r) 
 \longrightarrow   \text{Bew}_{\kappa}[\text{Neg}(Sb(r^{17}_{Z(x)})]}.
\end{equation}
\end{proof}

By (\ref{anti1}) and (\ref{anti2}), (\ref{15}) and (\ref{16}) turn out 
to be
false for any $x$ such that $x$ is a $PROOF$ $ARRAY$  which last FORMULA is $17\, \text{Gen}  \, r$ in $\kappa$, 
and the demonstration of Theorem VI 
cannot be accomplished.
Indeed, for $n$ such that $Bw_{\kappa}(n)$ and  $[l(n)]\, Gl\, n = 17\, \text{Gen}  \, 
r$, (\ref{anti1}) and (\ref{anti2}) 
yield in $\kappa$
\[
 \overline{ \overline{ n B_{\kappa}(17\, \text{Gen} \, r)} 
 \longrightarrow   \text{Bew}_{\kappa}[Sb(r^{17}_{Z(n)})]}
\]
\noindent and
\[
\overline{  n B_{\kappa}(17\, \text{Gen} \, r) 
 \longrightarrow   \text{Bew}_{\kappa}[\text{Neg}(Sb(r^{17}_{Z(n)})]},
\]
hence, within the case 

$\, $ ``\emph{1. $17\, \text{Gen} \, r$ is not $\kappa \text{--PROVABLE}$"} 
$\, $

\noindent in G\"{o}del's 1931 \cite[176-177]{godel1}, $ \text{Bew}_{\kappa}[\text{Neg}(Sb(r^{17}_{Z(n)})]$ 
has no basis to be obtained from 
$  n B_{\kappa}(17\, \text{Gen} \, r) $ in $\kappa$, therefore no 
proof that 
$\; 17\, \text{Gen} \, r$ is not $\kappa \text{--PROVABLE}\;$ can be 
achieved.

Accordingly, $(n)\,\overline{  n B_{\kappa}(17\, \text{Gen} \, r)}$ has no 
soundness  as a consequence of (6.1)  within the next case

$\, $ ``\emph{2.} $\text{Neg}(17\, \text{Gen} \, r)$ \emph{ is not $\kappa 
\text{--PROVABLE}$}" $\, $

\noindent \cite[176-177]{godel1}.
 Moreover,  for $n$ such 
that $Bw_{\kappa}(n)$  and $[l(n)]\, Gl\, n = 17\, \text{Gen}  \, r$, 
$\text{Bew}_{\kappa}[(Sb(r^{17}_{Z(n)})]$ 
is not a consequence of $\overline{  n B_{\kappa}(17\, \text{Gen} \, 
r)}$ in $\kappa$,  so that  neither a demonstration that 
$\; \text{Neg}(17\, \text{Gen} \, r)$ is not $\kappa 
\text{--PROVABLE} \; $
can  be  accomplished.

Consequently the statement of Theorem VI, ``\emph{For every} \emph{ 
$\omega$--consistent recursive class $\kappa$ of $FORMULAS$ there are 
recursive $CLASS$ $SIGNS$ $r$ such that neither $v \, \text{Gen} \, r$ 
nor $\text{Neg}(v \, \text{Gen} \, r)$ belongs to $Flg(\kappa )$ (where $v$ is the 
$FREE$ $VARIABLE$ of $r$})" \cite[172-173]{godel1}, has no proof.

Furthermore, the assertion that ``\emph{ it suffices for the existence of 
propositions undecidable that the class 
$\kappa$ be $\omega$--consistent}", is now meaningless \cite[176-177]{godel1}.
We can then state the following theorem.

\bigskip 

\begin{theorem}\label{no-undec}   
    The existence of undecidable propositions of the form $v \, \text{Gen}  
    \,    r$ is not a theorem in $\kappa$.
 \end{theorem}  
 
 \bigskip 
 
  The invalidity of theorems 
VIII, IX and XI, all consequent of theorem VI, follows immediately \cite[184-194]{godel1}.
In the outlined derivation of theorem XI \cite[192-194]{godel1}, 
the assertion that $17\, \text{Gen} \, 
r$  is not $\kappa \text{--PROVABLE}$, together with the assumption 
about the consistency of $\kappa$, have now no justification. 
Therefore the statements
$  
\text{Wid}_{s}(\kappa) \longrightarrow  \text{Bew}_{\kappa}(17\, \text{Gen} \, r)  
$ (23)
and
$  
\text{Wid}_{s}(\kappa) \longrightarrow  (x)Q(x,p) 
$ (24)
are in $P$ without soundness ($\text{Wid}_{s}(\kappa)¥$ means ``$\kappa$ 
is consistent" and is defined by $\text{Wid}_{s}(\kappa)$ $ \equiv $ $(Ex)(Form(x) \& 
\overline{\text{Bew}_{\kappa}¥}(x)$), see note 63 \cite[192-193]{godel1}).
\section{}\label{finalsynt}
To resume briefly, in G\"{o}del's 1931,   (15) and (16) 
are derived by means of definition (8.1) $Q(x,y)$ from (3) 
and (4) taking respectively $Q(x,y)$ as an instance of $R(x_{1},\ldots , x_{n})$ and 
$\overline{Q}(x,y)$ as an instance of $\overline{R}(x_{1},\ldots , x_{n})$, 
and  then the two cases ``\emph{1. $17\, \text{Gen} \, r$ is not $\kappa 
\text{--PROVABLE}$}" and ``\emph{2.} $\text{Neg}(17\, \text{Gen} \, r)$ 
\emph{ is} \emph{ not $\kappa 
\text{--PROVABLE}$}" 
are  derived respectively from (15) and (16) (figure \ref{fig1}).

Any deduction is sound at whatever row  if it is sound at all the preceding rows. 
Within the whole of G\"{o}del's 1931 the propositions (15) and (16) 
 are  meant to be sound to secure the 
correctness of the conclusions.
As  showed in Sections  \ref{wider} and 
 \ref{antig} this is not the case. 
By Lemmas (6) and (7), given  $Q(x,y)$, the 
negation of (15), (I),
can be deduced
 and, given $\overline{Q}(x,y)$, we can obtain the negation of (16), (II).
This invalidates the deduction of propositions (15) and (16) (figure 
\ref{fig2}). 
Indeed the soundness of 
``\emph{1. $17\, \text{Gen} \, r$ is not $\kappa \text{--PROVABLE}$}" and 
``\emph{2. } $\text{Neg}(17\, \text{Gen} \, r)$ \emph{ is} \emph{not $\kappa \text{--PROVABLE}$}"
is based respectively on the soundness of (15) and (16).
As we showed (15) and (16) are false (figure 
\ref{fig3}), therefore 
 there are no sound  bases in G\"{o}del's 1931 to state Theorem VI, Theorem (\ref{no-undec}). 

 A definition of a new relation, like 8.1, is supposed to be true at 
 the row it is asserted. Let us consider this definition at the light 
 of our results. For $Q(x,y)$ as a true relation we have both
\[ Q(x,y) \;\; \longrightarrow  \;\; \overline{  x B_{\kappa}(17\, \text{Gen} \, r)}  
 \longrightarrow   \text{Bew}_{\kappa}[Sb(r^{17}_{Z(x)})],\] 
\[ Q(x,y)  \;\; \longrightarrow  \;\; \overline{ \overline{  x B_{\kappa}(17\, \text{Gen} \, r)} 
 \longrightarrow   \text{Bew}_{\kappa}[Sb(r^{17}_{Z(x)})]}, \] 
 i. e. definition (8.1) implies a formula and its negation.
 A ancient tradition identifies a contradiction with a couple of 
 propositions where one is the negation of the other. If $Q(x,y)$ is 
 considered to be true then a contradiction follows, the couple (15) 
 and (I). Such a couple can be regarded from two distinguished point of 
 views: the 
 consequences for $Q(x,y)$ and the consequences for $P$ (and $\kappa$).  The
 first point-view delivers directly to the antinomic features of $Q(x,y)$. 
 The couple (15) and (I) is always false and can be derived when $Q(x,y)$ is
 supposed to be true, i.  e. $Q(x,y) \rightarrow \mathit{ false }$.  As a
 first immediate consequence, $Q(x,y)$ is false too.  But as showed above if
 $\overline{Q}(x,y)$ then we have both (16) and (II), a contradictory couple
 again.  All that clarify the antinomic features of 8.1 and leads to
 re-consider it in terms of the theory of definition along the lines already
 developed in \cite{catta3}.  From the second point-view of the consequences
 for $P$, if $Q(x,y) \rightarrow \mathit{ false }$ then {\it ex falso
 sequitur quodlibet}, $(\text{\tiny{$\sim$}}a \rightarrow (a \rightarrow b))
 $.  The existence of the undecidable formula in G\"{o}del's 1931 follows
 {\it ex falso}, i.  e. it follows because everything follows from false. 
 The law $(\text{\tiny{$\sim$}}a \rightarrow (a \rightarrow b)) $ is
 intuitionistically accepted but it leads $P$ to collapse.  $P$ is
 inconsistent, in it everything follows, even the existence of a property for
 which \emph{it is possible neither to give a counterexample nor to prove that it
 holds of all numbers}.

 \bigskip
 
 \bigskip

  \bibliographystyle{plain}


\bigskip

\begin{figure}[ht]
$
\begin{array}{ll}
(1) &  \dots  \\
(2) &   \dots     \\
(3) &    R(x_{1},\ldots , x_{n}) \longrightarrow   
\text{Bew}_{\kappa}[Sb(r^{u_{1}\ldots u_{n}}_{Z(x_{1})\ldots Z(x_{n})})]  \\
(4) &    \overline{R}(x_{1},\ldots , x_{n})  \longrightarrow   
\text{Bew}_{\kappa}[\text{Neg}(Sb(r^{u_{1}\ldots u_{n}}_{Z(x_{1})\ldots Z(x_{n})})]   \\
(5) &   \dots     \\
(6) &    \dots    \\
(6.1) &   \dots   \\
(7) &     \dots   \\
(8) &     \dots   \\
(8.1) & Q(x,y)     \\
& Q(x,y)  \longrightarrow \Big[\;\overline{  x B_{\kappa}[Sb(y^{19}_{ Z(y)})]} 
 \longrightarrow   \text{Bew}_{\kappa}[ Sb(q^{17 \quad 19}_{Z(x) 
 Z(y)})]\Big]\;\\
& \overline{Q}(x,y)  \longrightarrow  \Big[\;  x B_{\kappa}[Sb(y^{19}_{ Z(y)})] 
 \longrightarrow   \text{Bew}_{\kappa}[\text{Neg}( Sb(q^{17 \quad 
 19}_{Z(x) Z(y)}))] \Big]\;   \\
(9) &   \overline{  x B_{\kappa}[Sb(y^{19}_{ Z(y)})]} 
 \longrightarrow   \text{Bew}_{\kappa}[ Sb(q^{17 \quad 19}_{Z(x) Z(y)})]  \\
(10) &   x B_{\kappa}[Sb(y^{19}_{ Z(y)})] 
 \longrightarrow   \text{Bew}_{\kappa}[\text{Neg}( Sb(q^{17 \quad 19}_{Z(x) Z(y)}))]    \\
(11) &    \dots    \\
(12) &     \dots    \\
& (11), (12)  \longrightarrow  (13) \\
&  (12)  \longrightarrow  (14) \\
(13) &   \dots      \\
(14) &    \dots    \\
&(13), (14), (9) \longrightarrow   
 \Big[\;\overline{ x B_{\kappa}(17\, \text{Gen} \, r)} 
 \longrightarrow   \text{Bew}_{\kappa}[Sb(r^{17}_{Z(x)})]\Big]\;\\
&(13), (14), (10) \longrightarrow   
 \Big[\; x B_{\kappa}(17\, \text{Gen} \, r) 
 \longrightarrow   \text{Bew}_{\kappa}[Neg(Sb(r^{17}_{Z(x)}))]\Big]\;\\
(15) &   \overline{ x B_{\kappa}(17\, \text{Gen} \, r)} 
 \longrightarrow   \text{Bew}_{\kappa}[Sb(r^{17}_{Z(x)})]  \\
(16) &  x B_{\kappa}(17\, \text{Gen} \, r) 
 \longrightarrow   \text{Bew}_{\kappa}[Neg(Sb(r^{17}_{Z(x)}))]   \\
  &(15)    \longrightarrow  ``\emph{1. $17\, \text{Gen} \, r$ is not $\kappa \text{--PROVABLE}$}" \\
  &(16)  \longrightarrow ``\emph{2. } $\text{Neg}(17\, \text{Gen} \, r)$ \emph{ is not $\kappa 
\text{--PROVABLE}$}" \\ 
\\
\hline \\
   &   ``\emph{1. $17\, \text{Gen} \, r$ is not $\kappa \text{--PROVABLE}$}"\\
   & ``\emph{2. } $\text{Neg}(17\, \text{Gen} \, r)$ \emph{ is not $\kappa 
\text{--PROVABLE}$}" \\
\end{array}
$
\caption{ G\"{o}del's 1931}
\label{fig1}
\end{figure}

\begin{figure}[ht]
$
\begin{array}{ll}
(1) &  \dots  \\
(2) &   \dots     \\
(3) &    \dots   \\
(4) &    \dots    \\
(5) &   \dots     \\
(6) &    \dots    \\
(6.1) &   \dots   \\
(7) &     \dots   \\
(8) &     \dots   \\
(8.1) & Q(x,y)     \\
& Q(x,y)  \longrightarrow \Big[\;\overline{  x B_{\kappa}[Sb(y^{19}_{ Z(y)})]} 
 \longrightarrow   \text{Bew}_{\kappa}[ Sb(q^{17 \quad 19}_{Z(x) 
 Z(y)})]\Big]\;\\
& \overline{Q}(x,y)  \longrightarrow  \Big[\;  x B_{\kappa}[Sb(y^{19}_{ Z(y)})] 
 \longrightarrow   \text{Bew}_{\kappa}[\text{Neg}( Sb(q^{17 \quad 
 19}_{Z(x) Z(y)}))]   \Big]\; \\
(9) &   \overline{  x B_{\kappa}[Sb(y^{19}_{ Z(y)})]} 
 \longrightarrow   \text{Bew}_{\kappa}[ Sb(q^{17 \quad 19}_{Z(x) Z(y)})]  \\
(10) &   x B_{\kappa}[Sb(y^{19}_{ Z(y)})] 
 \longrightarrow   \text{Bew}_{\kappa}[\text{Neg}( Sb(q^{17 \quad 19}_{Z(x) Z(y)}))]    \\
(11) &    \dots    \\
(12) &    \dots     \\
& (11), (12)  \longrightarrow  (13) \\
&  (12)  \longrightarrow  (14) \\
(13) &    \dots     \\
(14) &    \dots    \\
&(13), (14), (9) \longrightarrow   
 \Big[\;\overline{ x B_{\kappa}(17\, \text{Gen} \, r)} 
 \longrightarrow   \text{Bew}_{\kappa}[Sb(r^{17}_{Z(x)})]\Big]\;\\
&(13), (14), (10) \longrightarrow   
 \Big[\; x B_{\kappa}(17\, \text{Gen} \, r) 
 \longrightarrow   \text{Bew}_{\kappa}[Neg(Sb(r^{17}_{Z(x)}))]\Big]\;\\
 \\
\hline \\
(15) &   \overline{ x B_{\kappa}(17\, \text{Gen} \, r)} 
 \longrightarrow   \text{Bew}_{\kappa}[Sb(r^{17}_{Z(x)})]  \\
(16) &  x B_{\kappa}(17\, \text{Gen} \, r) 
 \longrightarrow   \text{Bew}_{\kappa}[Neg(Sb(r^{17}_{Z(x)}))]   \\
\end{array}
$
\caption{Derivation of (15) and (16) within  G\"{o}del's 1931}
\label{fig2}
\end{figure}
\begin{figure}[ht]
$
\begin{array}{ll}
(1) &  \dots  \\
(2) &   \dots   \\
(3) &   \dots  \\
(4) &   \dots   \\
(5) &   \dots   \\
(6) &   \dots   \\
(6.1) &  \dots  \\
(6.2) &  \dots  \\
(6.2.1) & \dots   \\
(7) &    \dots  \\
(7.1) &  \dots  \\
&\text{ Lemma } (6)\\
&\text{ Lemma } (7)\\
(8) & \dots       \\
(8.1) & Q(x,y)     \\
(11) &    \dots  \\
(12) &    \dots  \\
(13) &    \dots  \\
& Q(x,y)  \longrightarrow \overline{ \overline{ x B_{\kappa}(17\, \text{Gen} \, r)} 
 \longrightarrow   \text{Bew}_{\kappa}[Sb(r^{17}_{Z(x)})]}\\
& \overline{Q}(x,y)  \longrightarrow   \overline{  x B_{\kappa}(17\, \text{Gen} \, r) 
 \longrightarrow   \text{Bew}_{\kappa}[\text{Neg}(Sb(r^{17}_{Z(x)})]} \\
(I) &  \overline{ \overline{ x B_{\kappa}(17\, \text{Gen} \, r)} 
 \longrightarrow   \text{Bew}_{\kappa}[Sb(r^{17}_{Z(x)})]}  \\
(II) & \overline{  x B_{\kappa}(17\, \text{Gen} \, r) 
 \longrightarrow   \text{Bew}_{\kappa}[\text{Neg}(Sb(r^{17}_{Z(x)})]}   \\
\\
\hline \\
 &(15) \text{ is false}  \\
&(16) \text{ is false}  \\
\end{array}
$
\caption{(15) and (16) are false}
\label{fig3}
\end{figure}
\end{document}